\documentclass[a4paper]{article}

\usepackage[margin=1in]{geometry}

\usepackage[T1]{fontenc}
\newcommand{\changefont}[3]{
\fontfamily{#1} \fontseries{#2} \fontshape{#3} \selectfont}

\changefont{ptm}{m}{n}

\usepackage{setspace} 

\usepackage{graphicx}    

\usepackage{amsfonts,amssymb,amsmath,amsgen,amsopn,amsbsy,theorem,graphicx,epsfig}
\usepackage{graphics}
\usepackage{mathrsfs}
\usepackage{color}

\newtheorem{definition}{Definition}[section]

\long\def\symbolfootnote[#1]#2{\begingroup%
\def\thefootnote{\fnsymbol{footnote}}\footnote[#1]{#2}\endgroup} 

\begin{document}

\begin{center}
\Large \textbf{A randomly determined unpredictable function}
\end{center}

\begin{center}
\normalsize \textbf{Marat Akhmet$^{a,}\symbolfootnote[1]{Corresponding Author Tel.: +90 312 210 5355,  Fax: +90 312 210 2972, E-mail: marat@metu.edu.tr}$, Mehmet Onur Fen$^b$, Ejaily Milad Alejaily$^c$} \\
\vspace{0.2cm}
\textit{\textbf{\footnotesize$^a$Department of Mathematics, Middle East Technical University, 06800 Ankara, Turkey}} \\
\textit{\textbf{\footnotesize$^b$Department of Mathematics, TED University, 06420 Ankara, Turkey}} \\
\textit{\textbf{\footnotesize$^c$College of Engineering Technology, Houn, Libya}} 
\vspace{0.1cm}
\end{center}

\vspace{0.3cm}

\begin{center}
\textbf{Abstract}
\end{center}

\noindent\ignorespaces
Recently, we have introduced unpredictable oscillations, which are in the basis of Poincar\'{e} chaos.   For theoretical analysis as well as for applications,  it is necessary to provide   constructive examples of unpredictable functions. We have already   provided such functions utilizing orbits of the  logistic map,  and  in the present paper we suggest   another  way  of  construction of  the  functions   by  applying    the Bernoulli random process.  A simulation for a randomly determined unpredictable function is  provided.

\vspace{0.2cm}
 
\noindent\ignorespaces \textbf{Keywords:} Unpredictable function, Unpredictable sequence, Bernoulli process, Poincar\'{e} chaos,  Symbolic dynamics

\vspace{0.6cm}

\section{Introduction and preliminaries}

The theory of oscillations extremely rely on functions, which can be either tabulated or formalized. The ones in the second category are based first of all on the functions which are trigonometric, polynomials, hyperbolic trigonometric and others. All of them have been tabulated in computer memories. Next ones are functions, which can be presented as finite or infinite sums of the former ones. They are evaluated by developing software programs and very helpful in applications. Other are oscillations produced as solutions of differential  equations. There exists, even, the large class in the qualitative theory of differential equations, oscillatory differential equations. The solutions are approved as oscillations by special type of criteria for the existence. In this study, we focus on  functions which are shaped through   qualitative conditions of definitions. They    make the core of the  research  area  in the theory of dynamical systems, issued  by A. Poincar\'{e}, G. Birkhoff, and others. These are periodic, quasi-periodic, almost periodic oscillations, recurrent and Poisson stable orbits. A special type of Poisson stable orbit called an unpredictable trajectory, which leads to Poincar\'{e} chaos in the quasi-minimal set, was introduced in the paper \cite{AkhmetUnpredictable}. 
Moreover, the papers \cite{AkhmetPoincare,AkhmetExistence,AkhmetUnpredSolnDE,AkhmetUnpredSolnDEDisE} were concerned with unpredictable solutions of various types of quasi-linear differential equations. In the present paper, we introduce a new  way  for  unpredictable functions construction  benefiting from the dynamics associated with the discrete distribution \cite{Cinlar}. We consider the process with a finite number of possible outcomes to generate an unpredictable sequence. The sequence is then used to construct a continuous unpredictable function. 	 
Thus, unpredictable oscillations appeared as solutions of linear or quasi-linear differential equations, i.e., as outputs of the systems, provided that there is an unpredictable input. The natural question how it is possible to choose the inputs being unpredictable arises. For this  reason in the previous papers \cite{AkhmetPoincare,AkhmetExistence,AkhmetUnpredSolnDE,AkhmetUnpredSolnDEDisE}, we introduced unpredictable functions built by applying orbits of the logistic map, which were verified to be unpredictable sequences. One can confirm that in this way we utilize several other discrete equations with dynamics topologically equivalent to  symbolic dynamics. This is why,  they  are in some sense the same  as those functions, which have been already determined in our research. For that reason, the task of  construction of  new  unpredictable oscillations  is undertaken in the present paper. We utilize the two principal issues for the solution of the problem. The  first one is  that the set of all orbits of the symbolic dynamics   coincides with all possible sequences of the symbols. Moreover,  realizations of the Bernoulli random process altogether are the set of sequences. Consequently,   constructing an orbit of a random process, we  obtain  an orbit of the symbolic dynamics and  simulate a part of the unpredictable sequence.   Thus, we obtain that a single iteration of the Bernoulli shift is the same as a trial for the Bernoulli  process.

The next definitions are concerned with  unpredictable sequences and functions.

\begin{definition} (\cite{AkhmetUnpredSolnDE}) \label{def1}
A bounded sequence $ \{ \nu_k \}$, $k \in \mathbb{Z} $, in $ \mathbb R^p $ is called unpredictable if there exist a positive number $ \varepsilon_{0} $ and sequences $ \{ \zeta_n \} $, $ \{ \eta_n \} $, $ n \in \mathbb N $, of positive integers both of which diverge to infinity such that $ \| \nu_{k + \zeta_n} - \nu_k \| \to 0 $ as $ n \to \infty$ for each $k$ in bounded intervals of integers and $ \| \nu_{\zeta_n + \eta_n} - \nu_{\eta_n} \| \geq \varepsilon_{0} $ for each $ n \in \mathbb N $.
\end{definition}
	
\begin{definition} (\cite{AkhmetPoincare}) \label{def2}
A uniformly continuous and bounded function $h : \mathbb R \to \mathbb R^p$ is unpredictable if there exist positive numbers $\varepsilon_{0}$, $\sigma$ and sequences $\left\{t_{n}\right\}$, $\left\{u_{n}\right\}$ both of which diverge to infinity such that $h(t+t_{n})\to h(t)$ as $n\to \infty$ uniformly on compact subsets of $\mathbb R$ and $\|h(t+t_{n})-h(t)\|\geq \varepsilon_{0}$ for each $t\in  \left[u_{n}-\sigma, u_{n}+\sigma \right]$ and $n\in \mathbb N$.
\end{definition}

Consider the  space  $\Sigma_m$  of bi-infinite sequences  $\ldots i_{-2}i_{-1} \textbf{.} i_{0}i_1 i_2  \ldots$  on finite number of complex  numbers $a_1,\ldots,a_m,$  with  the metric 
\begin{equation} \label{DistMetric}
d(I,J) =  \sum_{k=-\infty}^{\infty} \frac{|i_k - j_k|}{2^{|k|}}, 
\end{equation}
where   $I = ( \ldots i_{-2}i_{-1} \textbf{.} i_{0}i_1 i_2  \ldots),  J = ( \ldots j_{-2}j_{-1} \textbf{.} j_{0}j_1 j_2  \ldots),$  and $\left|\centerdot\right|$  is the  absolute value. 
Introduce the Bernoulli shift $ \varphi : \Sigma_m \to \Sigma_m $ such that
\begin{equation} \label{MapDefn}
\varphi ((\ldots i_{-2}i_{-1} \textbf{.} i_{0}i_1 i_2   \ldots )) = (\ldots i_{-2}i_{-1}i_{0} \textbf{.} i_1i_2 i_3   \ldots ).
\end{equation}
The map $\varphi$ is continuous and the metric space $\Sigma_m$ is compact \cite{Wiggins88}.

Let  us, now,  build an  unpredictable point   for  the dynamics   $ (\Sigma_m, d, \varphi) $ . Without loss of generality, we consider a particular case of the space when $m=2,$ $a_1 = 0,$ $a_2= 1.$  We need a collection of finite sequences $i^{r}_{k},$ $r \in \mathbb N,$ $k=1,2,\ldots, 2^r,$ consisting of $0$'s and $1$'s. Let us use the notations $i^1_1=(0)$ and $i^1_2=(1)$ for the sequences of length $1$. For each natural number $r$, we recursively define $i^{r+1}_{2k-1}=(i_k^r 0)$ and $i^{r+1}_{2k}=(i_k^r 1),$ $k=1,2,\ldots,2^r$, where $i^{r+1}_{2k-1}$ and $i^{r+1}_{2k}$ are obtained by respectively inserting $0$ and $1$ to the end of the sequence $i^r_k$ of length $r$.
For instance, $i^2_1=(i^1_1 0)=(00)$, $i^2_2=(i^1_1 1)=(01)$, $i^2_3=(i^1_2 0)=(10)$, and $i^2_4=(i^1_2 1)=(11)$ are the sequences of length $2$.    Now consider the following sequence  $
i^* = (\ldots i_8^3i_6^3i_4^3i_2^3i_4^2i_2^2.i_1^1i_1^2i_3^2i_1^3i_3^3i_5^3i_7^3\ldots) .$ In \cite{AkhmetUnpredictable} it  was proved that   $i^*$ is an unpredictable point of the dynamics.  

Because the trajectory  which initiates at   $i^*$ is dense in the quasi-minimal set $\Sigma_m$, the dynamics    is Poincar\'{e} chaotic according to Theorem 3.1 presented in paper \cite{AkhmetUnpredictable}. Moreover, there  is an  uncountable set  of  unpredictable points  in the set. From this discussion it implies that any numerical simulation of a discrete finite distribution is an approximation of an unpredictable sequence.  Indeed, the metric peculiarity  implies that  if one considers   the   point $i^*$   in $\Sigma_m$   as a  bi-infinite  sequence,  then it  is  easily  seen that  it  is an  unpredictable sequence in the  sense of Definition \ref{def1}.     This is in the base of the construction of an unpredictable function in the next section.

\section{Main result}

Let  us fix    a finite  string  $i_k,\ldots,i_p,  1\le k < p,$  on the  set  of    complex  numbers  $a_1,\ldots,a_m.$     It  can  be   accepted as an  arc of   a sequence from $\Sigma_m.$ Since of the  last  section discussion,  the string  can  be approximated  with  arbitrary  precision by  a shift  of the sequence $i^*.$     This  possibility  to  approximate  by  shifts  of the orbits  is  the  main advantage of the Poincar\'{e} chaos against  other  types of chaos. Taking  into  account  that  there   are  limits  for the approximation validity  in numerical   simulations by  computers, we can  admit  that  simulation of the  string  is simulation of the unpredictable sequence itself.   This  is why,  we accept that  finite   realizations of the Bernoulli  process,  which  are obtained randomly present   the unpredictable  sequence,  since,  at  first,   they  are not  periodic  even on a sufficiently  large interval  of  discrete  time,  and, secondly,  since of the above explanation the simulation is an  approximation of the  sequence  with arbitrary  precision.   The arbitrariness  guarantees    that  in applications we can  get  the simulations as the  unpredictable  sequence with  the atributes  listed in the definition.     Moreover, we must  not  be confused with the approximations in the  basis of the definition.  This  is true for  all types of functions, which  are determined  through infinitely  long   algorithms such as series,  for  instance. 

Fix   an   unpredictable sequence $i^*,$   which  is defined on  the  two real   numbers $a$   and $b.$   One can  find that  the unpredictability constant  $\epsilon_0$  can be taken  equal  to  $|a-b|.$   Define the function $\chi(t): \mathbb R \to \mathbb R$ through the equation
\begin{equation}  \label{func_unp}
\chi(t)=\displaystyle \int_{-\infty}^{t} e^{-(t-s)} \pi(s) ds,
\end{equation} 
where $\pi(t):\mathbb R \to \mathbb R$ is the piecewise constant function satisfying $\pi(t) = i^*_k$ for $t\in [k,k+1)$, $k \in \mathbb Z $. One can confirm that  $\displaystyle \sup_{t\in\mathbb R} \left|\chi(t)\right| \leq M_{\chi}$, where $M_{\chi}=\max\{|a|,|b|\}$.  

We will show that the function $\chi(t)$ defined by (\ref{func_unp}) is unpredictable. Consider a fixed compact interval $[\alpha, \beta]$ and a positive number $\varepsilon$. We assume without loss of generality that $\alpha$ and $ \beta$ are integers. Let us fix a positive number $\xi$ and an integer $\gamma <\alpha$ which satisfy the inequalities $Me^{-2(\alpha-\gamma)}< \varepsilon/4$ and $\xi(1-e^{-2(\beta-\gamma)})< \varepsilon.$
Suppose that $n$ is a sufficiently large natural number satisfying $\left|\pi(t+\zeta_n) - \pi(t)\right|< \displaystyle \xi$ for every $t$ in $[\gamma, \beta].$ Accordingly, we have for $t\in [\alpha,\beta]$ that
\begin{eqnarray*}
	|\chi(t+\zeta_n)-\chi(t)| & = &  \le  \displaystyle \int_{-\infty}^{\gamma} e^{-2(t-s)} \left|\pi(s+\zeta_n)-\pi(s)\right| ds + \int_{\gamma}^{\beta} e^{-2(t-s)}  \left|\pi(s+\zeta_n)-\pi(s) \right|  ds  \\
	&& \le \int_{-\infty}^{\gamma}e^{-2(t-s)}2ds+\int_{\gamma}^{\beta}e^{-2(t-s)}\xi ds \\
	&& \le 2Me^{-2(\alpha-\gamma)}+\frac{\xi}{2}[1-e^{-2(\beta-\gamma)}] \\
	&& < \varepsilon.
\end{eqnarray*}

Thus, 
$\left|\chi(t+\zeta_n) -  \chi(t) \right|  \to 0$ as $n \to \infty$
uniformly on the interval $[\alpha, \beta].$

Let us fix a number $ n$ and consider two alternative cases: (i) $|\chi(\eta_n+\zeta_n)-\chi(\eta_n)|< \frac{\epsilon_0}{8}$ and (ii) $|\chi(\eta_n+\zeta_n)-\chi(\eta_n)|\geq \frac{\epsilon_0}{8}.$

(i) There exists a positive number $\kappa<1$ such that
$e^{-2\kappa} = \frac{2}{3}.$ 
Using the relation
\begin{eqnarray}\label{theta}
&& \chi(t+\zeta_n)-\chi(t) = \chi(\eta_n+\zeta_n) - \chi(\eta_n) + \displaystyle   \int_{\eta_n}^t e^{-2(t-s)} (\pi(s+\zeta_n) -\pi(s)) ds  
\end{eqnarray}
we obtain that
\begin{eqnarray*}
	&&|\chi(t+\zeta_n)-\chi(t)|\geq |\displaystyle   \int_{\eta_n}^t e^{-2(t-s)} |\pi(s+\zeta_n) -\pi(s)) |ds-|\pi(\eta_n+\zeta_n) - \pi(\eta_n)|\geq\\
	&& \int_{\eta_n}^t e^{-2(t-s)}\epsilon_0 ds- \frac{\epsilon_0}{8} \geq \frac{\epsilon_0}{2}(1-  e^{-2\kappa})-\frac{\epsilon_0}{8} = \frac{\epsilon_0}{24}
\end{eqnarray*}
for $t \in [\eta_n + \kappa, \eta_n+1).$

(ii)  There exists a positive number $\kappa<1$ such that
$ 1-e^{-2\kappa} = \frac{\epsilon_0}{12}.$   From the relation (\ref{theta}) we get
\begin{eqnarray*}
	&&|\chi(t+\zeta_n)-\chi(t)|\geq |\chi(\eta_n+\zeta_n) - \chi(\eta_n)|-|\displaystyle   \int_{\eta_n}^t e^{-2(t-s)} (\pi(s+\zeta_n) -\pi(s)) ds| \geq \\
	&& \frac{\epsilon_0}{8} -\int_{\eta_n}^t e^{-2(t-s)}2ds \geq \frac{\epsilon_0}{8}-[1-e^{-2\kappa}]  = \frac{\epsilon_0}{24}	
\end{eqnarray*}
for $t \in [\eta_n, \eta_n+\kappa).$

Thus, $\chi(t)$ is  an  unpredictable   function.

It  is easy  to  see that $\chi(t)$ is a solution of the differential equation
\begin{eqnarray} \label{difequatchi}
x'=-x+h(t)
\end{eqnarray}	
with 
\begin{eqnarray} \label{initialchi}
\chi(0)=\int_{-\infty}^0 e^sh(s) ds,
\end{eqnarray}		
but we do not know the value $\chi(0)$ precisely, since it cannot be evaluated by the improper integral (\ref{initialchi}). Nevertheless, we utilize that $\chi(t)$ is an exponentially stable solution of 	equation (\ref{difequatchi}). Therefore, any solution $\varphi(t)$ of (\ref{difequatchi}) approximates $\chi(t)$. The approximation is better for larger $t$ such that $\left\|\chi(t)-\varphi(t)\right\|\leq \left\|\chi(0)-\varphi(0)\right\|e^{-t}$, $t \geq 0$. For that reason we take $\varphi(0)=0.5$ so that $\left\|\chi(t)-\varphi(t)\right\|\leq e^{-50}<10^{-17}$ for $t \in [50,100]$. It is less than Matlab precision between $50$ and $100$. Hence, the part of the time series of $\varphi(t)$ for $50 \leq t \leq 100$ can be accepted as the graph of the function $\chi(t)$.

\begin{figure}[ht!]
	\centering
	\includegraphics[width=1.0\linewidth]{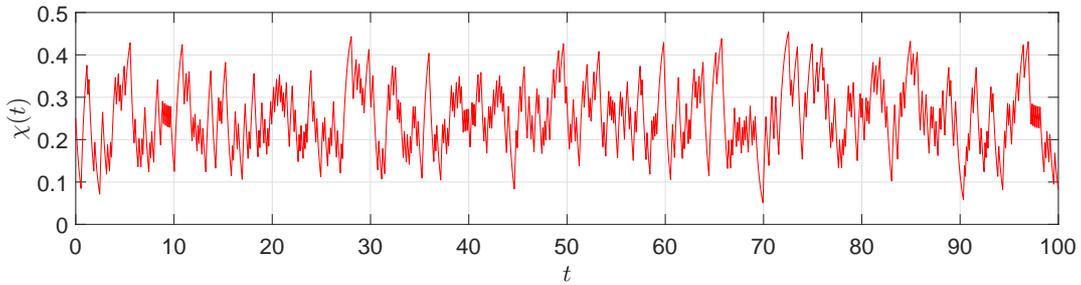}
	\caption{Time series of the unpredictable function $\chi(t)$}
	\label{UnpFunc}	
\end{figure}

In Figure \ref{UnpFunc} we depict the unpredictable function $ \chi(t) $ defined by equation (\ref{func_unp}). For the simulation, we use the function $ \pi(t) = i_k, \;  t\in [\mu(k-1),\mu k), \, \mu=0.1, \, k \in \mathbb N $. The sequence $ i_k $ is generated randomly such that $ i_k = 0, 1$ for each $k=1, 2, \ldots $.

\end{document}